\documentclass[letterpaper, 10 pt, conference]{ieeeconf}  

\IEEEoverridecommandlockouts                              

\usepackage{graphicx}
\usepackage{tikz}
\usepackage{caption}
\graphicspath{
    {figures/}
}

\usepackage{amsmath,amssymb,enumerate}

\usepackage{amsthm}
\usepackage{setspace}
\usepackage{booktabs}
\usepackage[usenames,dvipsnames,svgnames,table]{xcolor}
\usepackage{mathtools}
\usepackage{algorithm, algorithmicx, algpseudocode}
\usepackage{blindtext}
\usepackage{gensymb}
\usepackage{xparse}
\usepackage{lipsum}
\usepackage{mathrsfs}
\usepackage[mathscr]{euscript}
\usepackage{times}
\usepackage{cite} 
\usepackage{multicol}
\definecolor{darkgreen}{rgb}{0,0.6,0}
\usepackage[bookmarks=true,colorlinks=true,pdfpagemode=UseNone,citecolor=darkgreen,linkcolor=black,urlcolor=BrickRed,pagebackref]{hyperref}
\usepackage[caption=false,font=footnotesize]{subfig}
\usepackage{amsfonts}
\usepackage[utf8]{inputenc}
\usepackage[T1]{fontenc}
\usepackage{textcomp}
\usepackage{arydshln}
\usepackage{balance}
\usepackage{float}

\newtheorem{problem}{Problem}

\renewcommand{\thefigure}{\arabic{figure}}

\definecolor{note}{rgb}{0.1,0.1,1}
\definecolor{rephase}{rgb}{0.15,0.7,0.15}
\definecolor{bag}{rgb}{0.6,0.6,0.2}

\makeatletter
\renewcommand*\env@matrix[1][c]{\hskip -\arraycolsep
  \let\@ifnextchar\new@ifnextchar
  \array{*\c@MaxMatrixCols #1}}
\makeatother

\makeatletter
\newcommand{\mathleft}{\@fleqntrue\@mathmargin0pt}
\newcommand{\mathcenter}{\@fleqnfalse}
\makeatother

\definecolor{orange}{RGB}{255,127,0}

\title{\LARGE \bf
Switching Network System Identification via Convex Optimizations
}

\author{Kaito Iwasaki, Anthony Bloch, Maani Ghaffari
\thanks{Kaito Iwasaki, Anthony Bloch, Maani Ghaffari are with the University of Michigan, Ann Arbor, MI 48109, USA. \texttt{\{kaitoi, abloch, maanigj\}@umich.edu}.
}
}

\begin{document}

\maketitle
\thispagestyle{empty}
\pagestyle{empty}


\begin{abstract}
This paper introduces a convex optimization framework for identifying switched network systems, in which both the node dynamics and the underlying graph topology switch between a finite number of configurations. Building on our recent convex identification method for general switching systems, we extend the formulation to structured network systems where each mode corresponds to a distinct adjacency matrix. We show that both the continuous node dynamics and binary network topologies can be identified from sampled state–velocity data by solving a sequence of convex programs. The proposed framework provides a unified and scalable way to recover piecewise network structures from data without a prior knowledge of mode labels at each state. Numerical results on diffusively coupled oscillators demonstrate accurate recovery of both mode dynamics and switching graphs.
\end{abstract}

\IEEEpeerreviewmaketitle


\section{Introduction and motivation}
This paper introduces a convex optimization-based method for identifying \emph{switching network systems}, a class of hybrid dynamical systems in which both the nodal dynamics and the network topology evolve through discrete mode transitions. The approach is based on our recent work on \cite{iwasaki2025learning}, and applies it to the network dynamical setting studied in \cite{mouyebe2025network}, where each mode corresponds to a distinct adjacency matrix governing the coupling between nodes.

Network dynamical systems appear in a wide range of domains, including biological and chemical networks, social and economic systems, and distributed robotics. In many classical models, the underlying graph topology is assumed to be fixed or to vary smoothly over time. However, real-world networks often undergo \emph{abrupt structural changes}, such as link activations, failures, or switching of communication patterns. Capturing such phenomena requires a hybrid framework that accounts for discrete transitions in both dynamics and connectivity.

The goal of this paper is therefore to develop a convex and data-driven formulation which, from observed trajectories of node dynamics, jointly identifies
\begin{enumerate}
    \item the continuous dynamics at each node and 
    \item the discrete switching structure of the network. 
\end{enumerate}
This formulation enables the recovery of piecewise network models without requiring pre-labeled mode information and provides a foundation for identification of switching network systems.

\section{Related work}
Network dynamical systems (NDS) describe the evolution of interacting agents induced by their coupling through a graph structure. Such frameworks model a wide range of real-world applications such as biological, engineering, and social networks such as cell synchronization, power grids, and consensus dynamics. There is a vast amount of literature on these subjects, and we cannot list them all, but some representative work may be found in \cite{Mirollo1990synchro, Acebron2005Kuramoto, Dorfler2013powergrids, Olfati2007consensus}. Moreover, the concept of graph-coupled dynamics extends beyond its classical setting to modern neural network architectures. For example, graph neural networks (GNNs) can be interpreted as discretized dynamical systems on graphs, where each layer performs a message-passing step analogous to the time evolution of coupled agents \cite{poli2020gnode}. 

In many practical scenarios, however, the underlying network topology or coupling parameters vary with time or operating conditions, giving rise to switching network dynamical systems (SNDS) that capture mode-dependent interactions and structural changes within the network. For example, \cite{moreau2005stability} gives graph-connectivity conditions for linear consensus under time-varying graph, and \cite{yang2015network} studies the synchronization and stability of NDS with time-dependent graph structures using Lyapunov-based synchronization results for nonlinear agents under switching interactions. In the  work \cite{mouyebe2025network}, the effects of the topology of the network on Lypaunov stability are investigated in both the static and switched setting, where in a switched network common Lyapunov functions are derived which can be used to prove stability. 

While these papers address stability and synchronization assuming a known network structure, many practical situations require identifying the underlying topology and coupling law directly from data. In the static setting, this has led to a large body of work on graph learning and network reconstruction, where one seeks to recover the adjacency or Laplacian matrix that best explains observed node signals. Convex formulations for learning graphs from smooth or diffusion-like signals have been developed in \cite{kalofolias2016learn, dong2019learn}, while reconstruction from dynamics data has been studied through regression-based approaches such as \cite{timme2014revealing, ching2015reconstruct}. These static-identification methods motivate us to extend such frameworks to switching networks, where both the topology and node dynamics can vary across modes. A recent paper investigates a switching system identification via convex optimization \cite{iwasaki2025learning}. Building on this framework, we extend these ideas to switching network dynamical systems, jointly recovering network topology and mode-dependent dynamics from data.
\section{Preliminaries}
\subsection{Network dynamical systems}
Let $\mathcal{G} = (V, E)$ be a graph of $N$ vertices in which each vertex $x_i \in V$ is a point in $\mathbb{R}^n$ governed by a combination of isolated dynamics 
\begin{equation}\label{eq: isolated dynamics}
\dot x_i = f(x_i)
\end{equation}
and some coupling interaction from adjacent vertices. Dynamical systems of this sort are referred to as \emph{network systems}. 

In this work, we consider a model of network systems of the form
\begin{equation} \label{eq: network dynamics}
    \dot x_i = f(x_i) + \sum_{j=1}^N \mathbf{A}_{ij}\phi(x_i, x_j), \quad 1 \le i \le N,
\end{equation}
where $\mathbf{A}$ is the adjacency matrix of the graph $\mathcal{G}$ referred to as the \emph{network topology}, and $\phi:\mathbb{R}^n \times \mathbb{R}^n \to \mathbb{R}^n$ is a coupling function that models the pairwise interaction of the vertices. In particular, we will consider \emph{linear couplings} which have the form
\begin{equation}\label{eq: linear coupling}
    \phi(u,v) = C(\alpha u + \beta v),
\end{equation}
where $C \in \mathbb{R}^{n \times n}$ is the \emph{coupling channel matrix} and $\alpha, \beta \in \mathbb{R}$ are coupling parameters. A common choice of such $C$ is a diagonal matrix with Boolean entries such as identity matrices. We can also write (\ref{eq: network dynamics}) using the Kronecker product as
\begin{equation}\label{eq: network dynamics coupled Laplacian}
    \dot{\mathbf{x}} = \mathbf{f}(\mathbf{x}) + (L_{\alpha\beta} \otimes C)\mathbf{x},
\end{equation}
where $\mathbf{x} = (x_1, \ldots, x_N)^T \in \mathbb{R}^{Nn}$ is the joint state, and $\mathbf{f}(\mathbf{x}) = (f(x_1), \ldots, f(x_N))^T$ is the joint vector field of the network system. The coupling Laplacian $L_{\alpha\beta}$ parametrized by $\alpha, \beta$ is defined as
\begin{equation}
    L_{\alpha\beta } = \alpha D^{in} + \beta \mathbf{A},
\end{equation}
where $D^{in} := \operatorname{diag}(\mathbf{A}\mathbf{1})$ is the \emph{in-degree matrix} of $\mathbf{A}$ where $\mathbf{1} = (1,\ldots,1)^T$ is a vector filled with ones. 

In the sequel, we make the following assumption: $\mathcal{G}$ contains no multiple edges between any two vertices. In other words, we restrict ourselves to the class of graphs whose adjacency matrices $\mathbf{A}$ have entries $0$ or $1$. Note that this assumption still allows a broad class of graphs, as they may not be simple graphs, i.e., we allow graphs to be undirected and/or have self-loops. 

\subsection{Switching systems}
A \emph{switching system} is a dynamical system whose vector field
switches among a finite set of modes according to a rule that may depend on
the state, time, or other internal/external variables. Formally, we write
\begin{equation}\label{eq:switching_general}
    \dot{x} = X_{\sigma(x,t)}(x), \quad \sigma(x,t) \in \mathcal{M} := \{1,2,\ldots,M\},
\end{equation}
where each mode $j \in \mathcal{M}$ is associated with a smooth vector field $X_j : \mathbb{R}^n \to \mathbb{R}^n$, and $\sigma(x,t)$ denotes the \emph{switching signal} determining which mode is active. The switching can be classified as \emph{state-dependent}, \emph{time-dependent}, or more generally \emph{event-triggered}. Throughout this paper we focus on state-dependent switching so that $\sigma(x,t)$ is independent of time.

For convenience, we represent the mode activation using binary indicators
\begin{equation}\label{eq: binary constraints}
    \lambda_j(x) \in \{0,1\}, \quad 
    \sum_{j=1}^{M}\lambda_j(x) = 1,
\end{equation}
so that \eqref{eq:switching_general} can be equivalently written as
\begin{equation}\label{eq:switching_sum}
    \dot{x} = \sum_{j=1}^{M}\lambda_j(x)X_j(x).
\end{equation}
To enable convex formulations, we write the binary constraints in \eqref{eq: binary constraints} as
\begin{equation}\label{eq:lambda_relax}
    \lambda_j(x)(\lambda_j(x) - 1) = 0, \quad 
    \sum_{j=1}^{M}\lambda_j(x) = 1.
\end{equation}
This later allows us to use moment relaxation to a mixed integer program that we formulate in Section similar to the framework introduced in \cite{iwasaki2025learning}.

\section{Switching network systems}
If the graphs model real-world systems such as biological networks or social network systems, the underlying network topology is not fixed but evolves dynamically in response to the system's internal state. For example, biological links can be activated or suppressed according to concentration levels, and social ties can form or dissolve depending on behavioral thresholds. These phenomena motivate us to study \emph{switching network dynamical systems}, where the network topology changes depending on the current state of the system. The stability analysis of such systems is studied in \cite{mouyebe2025network}. Here   we present such network systems as a special class of switching systems.
\subsection{Model description}
Suppose we have a finite collection of graph configurations $\{\mathcal{G}_1, \ldots, \mathcal{G}_M\}$ on the same vertex set $V$, where each graph $\mathcal{G}_j$ is governed by a linearly coupled network vector field $F_j: \mathbb{R}^{Nn} \to \mathbb{R}^{Nn}$ as in (\ref{eq: network dynamics coupled Laplacian})
\begin{equation}
    F_j(\mathbf{x}) = \mathbf{f}_j(\mathbf{x}) + (L_{\alpha\beta}^j \otimes C)\mathbf{x}.
\end{equation}
Note that the coupling parameters $\alpha,\beta$, and the coupling channel matrix $C$ are fixed. Then, the \emph{switching network dynamical systems} can be modeled as
\begin{equation}\label{eq: switching network systems}
    \dot{\mathbf{x}}(t) = F_{\sigma(\mathbf{x}(t))}(\mathbf{x}(t)),
\end{equation}
with the switching signal $\sigma: \mathbb{R}^{Nn} \to \mathcal{M}$ for the finite index set $\mathcal{M}$ of graphs. As in \eqref{eq:switching_sum}, we rewrite \eqref{eq: switching network systems} by taking the convex combination of mode dynamics
\begin{equation}\label{eq: switched network system cvx}
\begin{aligned}
    \dot{\mathbf{x}}(t) &= \sum_{j=1}^M \lambda_j(\mathbf{x}(t))F_j(\mathbf{x}(t)),\\
    \text{s.t.} \quad &\lambda_j(\mathbf{x}) \in \{0,1\}, \quad \sum_{j=1}^M \lambda_j(\mathbf{x}) = 1.
\end{aligned}
\end{equation}
As before, we rewrite the mixed-integer constraint as a polynomial constraint
\begin{equation}\label{eq: poly eqality constraints}
    \lambda_j(\mathbf{x})(\lambda_j(\mathbf{x})-1) = 1, \quad \sum_{j=1}^M \lambda_j(\mathbf{x}) = 1.
\end{equation}
\section{Switching network systems identification}
In this section, we introduce the method of identifying a switching network system and its graph configuration via convex optimizations. We assume that the coupling structure such as \eqref{eq: linear coupling} is known in advance.
\subsection{Identification objectives} 
Given sampled data $\{ \mathbf{x}_i, \dot{\mathbf{x}}_i \}_{i=1}^{N'}$, the goal of identification is to jointly recover:
\begin{enumerate}
    \item the \emph{mode-wise node dynamics} $\mathbf{f}_j(\mathbf{x})$,
    \item the \emph{mode-dependent graph topologies} $\{\mathbf{A}_j\}$,
    \item the \emph{switching rule} encoded by $\lambda_j(\mathbf{x})$.
\end{enumerate}
To formulate this as a tractable optimization problem, we treat the vectorized adjacency matrices
\begin{equation}\label{eq:adj_vec}
    a_j = \operatorname{vec}(\mathbf{A}_j) \in \{0,1\}^{(Nn)^2}
\end{equation}
as decision variables.  
Relaxing the Boolean constraint $a_j \in \{0,1\}^{(Nn)^2}$ to its convex hull, i.e., a linear simplex constraint
\begin{equation}\label{eq: simplex adj}
    a_j \in \Delta_a = [0,1]^{(Nn)^2}
\end{equation}
enables convex optimization formulations. The present framework leverages a combination of \emph{moment-based semidefinite relaxation} and \emph{linear simplex relaxation} as in \cite{iwasaki2025learning} to handle this mode-dependent network system. The details of this convex formulation are introduced in the next section.

After solving the relaxed program, the recovered adjacency vectors $a_j$ are projected back to the binary set $\{0,1\}^{(Nn)^2}$ by nearest-integer rounding, yielding estimated graph topologies $\hat{\mathbf{A}}_j$. At convergence, the procedure provides consistent estimates of (i) the network configurations $\{\hat{\mathcal{G}}_j\}$, (ii) the mode-specific node dynamics $\hat{\mathbf{f}}_j$, and (iii) the switching rule $\hat{\lambda}_j(\mathbf{x})$.

\subsection{Problem formulation}
Consider a set of measurements $\{\mathbf{x}^i(t), \dot{\mathbf{x}}^i(t)\}_{i=1}^{N'}$ of sampled joint states and velocities of switching network systems. Our goal is to identify the hidden switching rule and the network system simultaneously that best matches the measured samples. Denote $\lambda_j^i := \lambda_j(\mathbf{x}^i)$, and we ask the following:
\begin{problem}\label{prob: mip network id}
    Consider a cost function $c(\mathbf{x},\dot{\mathbf{x}}) \ge 0$ and require that $c(\mathbf{x}^i,\dot{\mathbf{x}}^i) = 0$ if $\dot{\mathbf{x}}^i = \sum_{j=1}^M \lambda_j^iF_j(\mathbf{x}^i)$ for true $F_j$ and $\lambda_j$. We solve the following minimization problem:
    \begin{equation}
        \min_{\lambda_j, f_j, \mathbf{A}_j} \quad \sum_{i=1}^{N'} c(\mathbf{x}^i, \dot{\mathbf{x}}^i).
    \end{equation}
    over all functions $f_j$ , all vectors $\lambda_j$ and all matrices $\mathbf{A}_j$ satisfying the constraints in \eqref{eq: switched network system cvx} and \eqref{eq:adj_vec}.
\end{problem}

We note that there are several issues for numerical implementation. First, searching $f_j$ over all possible functions is an  infinite-dimensional problem. Moreover, constraints on $\lambda_j$ and $\mathbf{A}_j$ makes the problem an integer programming problem which are computationally expensive to solve.

To overcome these difficulties, we reduce it  to a finite-dimensional setting by modeling $f_j$ as a polynomial approximation. Moreover, we use the method introduced in \cite{iwasaki2025learning} to relax the integer constraints on $\lambda_j$ and $\mathbf{A}_j$ to convex constraints. Consider the polynomial basis of an $n$ dimensional system $x \in \mathbb{R}^n$ of order up to $d$:
\begin{equation}\label{eq: monomial basis}
    \phi_d(x) = \left[1, x_1, x_2, \ldots, x_1^d, x_2^d,\ldots , x_n^d\right].
\end{equation}
We model the isolated dynamics as a polynomial function
\begin{equation}
    \mathbf{f}_j(x) = C_j \phi_d(x), \quad C_j \in \mathbb{R}^{n\times P},
\end{equation}
where $P$ is the length of $\phi_d(x)$ and $C_j$ is a matrix of polynomial coefficients. We first relax the integer constraints in \eqref{eq:adj_vec} to the simplex constraint in \eqref{eq: simplex adj}. Then, we use the moment-based semidefinite relaxation to the mixed-integer mode constraints \eqref{eq: switched network system cvx} instead of imposing the polynomial equality constraints in \eqref{eq: poly eqality constraints}. Specifically, we introduce a truncated sequence of moments 
$y = \{y_\alpha\}_{|\alpha|\le 2k}$, with the truncation bound $k \in \mathbb{N}$, corresponding to a measure supported on the feasible set
\begin{equation}\label{eq:simplex_support}
    G_\lambda = 
    \{\lambda \in \mathbb{R}^M : 
    \lambda_j(\lambda_j-1) = 0, \ \forall j , \ \mathbf{1}^T\lambda = 1\}.
\end{equation}
Let $M_r(y)$ denote the order-$r$ moment matrix, and for a constraint polynomial $p(\lambda)$, let $M_{r-d_p}(py)$ denote the associated \emph{localizing matrix}, where $d_p=\lceil\operatorname{deg} p/2\rceil$. In the current setting, the semidefinite relaxation of polynomial constraint \eqref{eq: poly eqality constraints} then takes the form
\begin{equation}\label{eq:moment_localizing_equalities}
\begin{aligned}
    & M_r(y) \succeq 0,\\
    & M_{r-1}(\lambda_j(\lambda_j-1)y) = 0, \quad j = 1,\ldots,M,\\
    & M_{r-1}(\left(\mathbf{1}^T \lambda-1\right) y)=0.
\end{aligned}    
\end{equation}
Following the flat-extension theory developed by Curto and Fialkow \cite{CurtoFialkow1996,CurtoFialkow1998}, low-order moment relaxations such as $r=1$ or $r=2$ often yield nearly atomic moment matrices, from which discrete mode assignments can be recovered via rank tests or spectral extraction, as later formalized in the polynomial optimization framework of Lasserre \cite{Lasserre2001Moments}. Readers interested in more theoretical background can refer to an introductory text \cite{Lasserre_2015} by Lasserre. More specifically, the details specialized to the current setting can be found in \cite{iwasaki2025learning}. For $r=1$, the moment conditions in \eqref{eq:moment_localizing_equalities} reduces to
\begin{equation}
\begin{bmatrix}
1 & \lambda^T \\
\lambda & \Lambda
\end{bmatrix} \succeq 0, \quad \Lambda_{j j}=\lambda_j, \quad \mathbf{1}^T \lambda=1
\end{equation}
where we used $\lambda$ to denote the first order moments $y_j$, and $\Lambda_{i j}=y_{i j}$.

To apply the moment- and linear simplex-based convex optimization scheme, we first reformulate the Problem \ref{prob: mip network id} as another mixed integer program with respect to the one-norm cost function of the residual of the fitted network system:
\begin{problem}
    We minimize the discrepancy between the measured time derivative and the fitted switching network system:
    \begin{equation}
        \begin{aligned}
        \min_{\lambda^i, C_j, a_j}  &\quad \sum_{i = 1}^{N'}  \left\|\dot{\mathbf{x}}_i - \sum_{j=1}^M \lambda_j^i \, F_j(\mathbf{x}^i) \right\|_1 \\
        \mathrm{s.t.}\quad &\ \lambda_j^i(1-\lambda_j^i) = 0, \quad \mathbf{1}^T\lambda^i = 1, \ \forall i, j, \\
        &\quad a_j \in \{0,1\}^{(Nn)^2},  \quad \forall j\\
        \quad& -\eta \le C_j \le \eta, \quad \forall j,
        \end{aligned}
    \end{equation}
    where the inequalities on $C_j$ are component-wise, and $\eta \in \mathbb{R}$ defines a box constraint to ensure the feasible set is compact. 
\end{problem}
Note that the $\ell_1$-cost in the above mixed integer program is nonconvex as it involves cross terms between decision variables $\lambda^i$, $C_j$, and $a_j$. However, when either $\lambda^i$ or $C_j$ and $a_j$ are fixed, the cost reduces to the one-norm of the linear function of the free variables. Hence, corresponding subproblems become convex programs by the above semidefinite and linear relaxations of mode and adjacency matrix constraints, and by the standard technique of introducing slack variables in the $\ell_1$-norm objective \cite{Teng2025ACC}. This results in an alternating convex optimization algorithm that optimize $\lambda^i$ and $a_j$ via a semidefinite program (SDP) and a linear program (LP).
\subsection{Mode identification}
Randomly initialize the entries $a_j = \bar{a}_j$ of the adjacency matrices and the polynomial coefficients $C_j = \bar{C}_j$. Introduce slack variables $\delta_k^i \geq 0$ for $k=1, \ldots, nN$. Let us denote $c_k^i(a,C,\lambda):= [\dot{\mathbf{x}}^i-\sum_{j=1}^M \lambda_j^i F_j(\mathbf{x}^i)]_k$. We then solve the following SDP using order-1 moment relaxation:
\begin{problem}(Mode search by SDP)\label{prob: mode search}
    Given fixed $\{\bar{a}_j,\bar{C}_j\}_{j=1}^M$, solve for $\{\lambda^i, \Lambda^i\}_{i=1}^{N'}$:
\begin{equation}\label{eq:mode_SDP_general}
\begin{aligned}
\min_{\lambda^i,\Lambda^i,\delta^i}\quad & \sum_{i=1}^{N'}\sum_{k=1}^{Nn} \delta_k^i \\
\text{s.t.}\quad
& -\delta_k^i \le c^i_k(\bar{a}, \bar{C},\lambda) \le \delta^i_k, \quad \delta^i_k \ge 0, \\
&\begin{bmatrix}
1 & (\lambda^i)^T \\
\lambda^i & \Lambda^i
\end{bmatrix} \succeq 0, \quad \operatorname{diag}(\Lambda^i) = \lambda^i,\quad \mathbf{1}^T \lambda^i = 1.
\end{aligned}
\end{equation}
\end{problem}
In practice, we post-process each identified mode assignment $\lambda^i$ by hardening to the standard basis vector $e_{j^\star} = [0, \ldots,1,\ldots,0]^T$ with $j^\star \in \arg\max_j \lambda_j^i$. Note that for relaxations of order $r > 1$, the inequality constraints for the slack variables will become LMIs for the corresponding order-$r$ localizing matrix constraints in \eqref{eq:moment_localizing_equalities}.
\subsection{Network dynamics identification}
After solving the Problem \ref{prob: mode search}, we fix obtained $\lambda^i = \bar{\lambda}^i$ and solve the following LP to search the network dynamics:
\begin{problem}[Network dynamics search by LP]\label{prob: dynamics search}
    Given fixed $\{\bar{\lambda}^i\}_{i=1}^{N'}$, solve for $\{a_j, C_j\}_{j=1}^M$:
    \begin{equation}\label{eq: network dynamics LP}
        \begin{aligned}
            \min_{a,C, \delta^i} \quad & \sum_{i=1}^{N'}\sum_{k=1}^{Nn} \delta^i_k \\
            \text{s.t} \quad & -\delta_k^i \le c_k^i(a,C,\bar{\lambda}) \le \delta_k^i, \quad \delta_k^i \ge 0,\\
            & a_j \in \Delta_a, \quad -\eta \le C_j \le \eta, \quad \forall j.
        \end{aligned}
    \end{equation}
\end{problem}
After solving for the optimal $a_j^\star$, we impose the hard integer constraint to update them to $\hat{a}_j$ to satisfy \eqref{eq:adj_vec}. That is, for instance for each entry $(\hat{a}_j)_m$ of $\hat{a}_j \in \mathbb{R}^{(Nn)^2}$, we set
\begin{equation}
    (\hat{a}_j)_m = \begin{cases}
        1, \quad \text{if } (a_j^\star)_m > 0.5,\\
        0, \quad \text{otherwise}.
    \end{cases}
\end{equation}
We alternate these subproblems until convergence to obtain the estimated $\hat{\lambda}^i$, $\hat{a}_j$, and $\hat{C}_j$.
\subsection{Graph reconstruction}
After convergence of the above alternating optimization, each binarized vector is then reshaped back into an adjacency matrix
\begin{equation}
    \hat{\mathbf{A}}_j = \operatorname{reshape}(\hat{a}_j, Nn, Nn),
\end{equation}
from which the corresponding identified graph $\hat{\mathcal{G}}_j = (V, E_j)$ is recovered by setting
\begin{equation}
    (p,q) \in E_j \quad \text{if and only if} \quad (\hat{\mathbf{A}}_j)_{pq} = 1.
\end{equation}
When the system is defined on $N$ nodes each of dimension $n$, we typically extract the block-diagonal structure of $\hat{\mathbf{A}}_j$ to obtain the node-to-node adjacency pattern. Small numerical noise around $0.5$ can be handled by adaptive thresholding or by enforcing sparsity penalties during identification. The resulting graphs $\{\hat{\mathcal{G}}_1,\ldots,\hat{\mathcal{G}}_M\}$ represent the estimated network topologies corresponding to each identified dynamical mode.
\section{Example and numerical analysis}
In this section we present a simple example of state-dependent switching network systems and a numerical simulation for the system identification. 

\begin{figure}[h]
    \centering
    \resizebox{0.9\linewidth}{!}{%
    \begin{tikzpicture}[
        node/.style={circle,draw,fill=blue!20,minimum size=0.45cm,inner sep=1pt},
        every edge/.style={very thick},
        bidir/.style={<->, thick},   
        dir/.style={->, thick}       
    ]
    \node[node](A1) at ( 0, 2) {\scriptsize 1};
    \node[node](B1) at (-1.6,-0.6) {\scriptsize 2};
    \node[node](C1) at ( 1.6,-0.6) {\scriptsize 3};

    \node at (0,3) {\(\mathcal{G}_1\): Complete graph \(K_{3}\)};

    \draw[bidir] (A1)--(B1);
    \draw[bidir] (A1)--(C1);
    \draw[bidir] (B1)--(C1);

    \node[node](A2) at (6, 2) {\scriptsize 1};
    \node[node](B2) at (4.4,-0.6) {\scriptsize 2};
    \node[node](C2) at (7.6,-0.6) {\scriptsize 3};

    \node at (6,3) {\(\mathcal{G}_2\): Directed \(3\)-cycle};

    \draw[dir] (A2) -- (B2);  
    \draw[dir] (B2) -- (C2);  
    \draw[dir] (C2) -- (A2);  

    \end{tikzpicture}
    }
    \caption{Illustrations of graph stricture of $K_3$ and $C_3$. Graphs switch between \(\mathcal{G}_1\) and \(\mathcal{G}_2\).}
    \label{fig: graphs}
\end{figure}
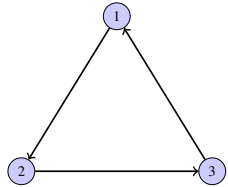

Consider a switching network dynamical system in $\mathbb{R}$ with 3 nodes and 2 modes. Let $\beta = - \alpha = -1$, i.e., the coupling is \emph{diffusive}. This gives the coupling Laplacian of the form 
\begin{equation}
    L_{\alpha\beta}^j = D_j^{in} - A_j.
\end{equation}
For simplicity, we will denote this Laplacian by $L_j$. Note that in some sense, this is the most classical form of the Laplacian coupling graph because the coupling $\phi(x_i, x_j) = x_j - x_i$ models  physical quantities as flowing from nodes with higher values to those with lower values in the most naive way, mirroring the Laplacian operator in heat/mass transport. Assume that the node dynamics is unchanged between the 2 graphs and are modeled by a single quadratic vector field $f(x) = \alpha_0 + \alpha_1 x + \alpha_2x^2$. Let the first graph $\mathcal{G}_1$ be a complete graph $K_3$, which is a fully connected simple undirected graph with no self-loops. Let the second graph $\mathcal{G}_2$ be a directed $3$-cycle $C_3$ where the edges point in one specific direction such as $1\to 2\to 3\to 1$. The switching happens when the joint state leaves/enters a prescribed compact set such as a ball of radius $R$. This can be written as
    \begin{equation}
        \mathcal{G}_j = \begin{cases}
            \mathcal{G}_1, & \text{if } x_1^2 + x_2^2 + x_3^2 \le R^2,\\
            \mathcal{G}_2, & \text{if } x_1^2 + x_2^2 + x_3^2 > R^2,
        \end{cases}
    \end{equation}
    and these graph structures are depicted in Figure \ref{fig: graphs}.

To derive its dynamics, we compute the coupling Laplacian for each graph. For $\mathcal{G}_1$, we have
    \[\begin{aligned}
        &\mathbf{A}_1 = \begin{bmatrix}
            0 & 1 & 1\\1 & 0 & 1\\ 1 & 1 & 0
        \end{bmatrix},\quad D_1^{in} = \operatorname{diag}(\mathbf{A}_1\mathbf{1}) = \begin{bmatrix}
            2 & 0 & 0\\ 0 & 2 & 0\\ 0 & 0 & 2
        \end{bmatrix}.
    \end{aligned}
    \]
    Thus, we have
    \[L_1 = D_1^{in} - \mathbf{A}_1 = \begin{bmatrix}
            2 & -1 & -1\\ -1 & 2 & -1\\ -1 & -1 & 2
        \end{bmatrix}.\]
    Similarly, for $\mathcal{G}_2$, we have
    \[\begin{aligned}
        &\mathbf{A}_2 = \begin{bmatrix}
            0 & 1 & 0\\0 & 0 & 1\\ 1 & 0 & 0
        \end{bmatrix},\quad D_2^{in} = \operatorname{diag}(\mathbf{A}_2\mathbf{1}) = \begin{bmatrix}
            1 & 0 & 0\\ 0 & 1 & 0\\ 0 & 0 & 1
        \end{bmatrix},
    \end{aligned}\]    
        so that
    \[L_2 = D_2^{in} - \mathbf{A}_2 = \begin{bmatrix}
        1 & -1 & 0\\ 0 & 1 & -1\\
        -1 & 0 & 1
    \end{bmatrix}.\]
    Let the joint state vector be $\mathbf{x}(t) = (x_1(t), x_2(t),x_3(t))^T$ and the joint node dynamics be denoted by $\mathbf{f}(\mathbf{x}) = (f(x_1),f(x_2), f(x_3))^T$. We set the coupling channel matrix (or scalar in this case) to $C = -1$. Then, the switching network system takes the form
    \begin{equation}
        \dot{\mathbf{x}}(t) = \mathbf{f}(\mathbf{x}) - \lambda_1(\mathbf{x}(t))L_1\mathbf{x}(t) - \lambda_2(\mathbf{x}(t))L_2 \mathbf{x}(t),
    \end{equation}
    such that
    \begin{equation}
        \lambda_1(\mathbf{x}(t)), \lambda_2(\mathbf{x}(t)) \in \{0, 1\}, \quad  \lambda_1(\mathbf{x}(t)) + \lambda_2(\mathbf{x}(t)) = 1 \quad\forall t,
    \end{equation}
    where $ \lambda_1(\mathbf{x}) = 1$ if $x_1^2 + x_2^2 + x_3^2 \le R^2$, and $0$ otherwise, and $\lambda_2(\mathbf{x})$ is defined similarly.
    
\begin{figure}[ht!]
      \centering
      \includegraphics[width=0.8\linewidth]{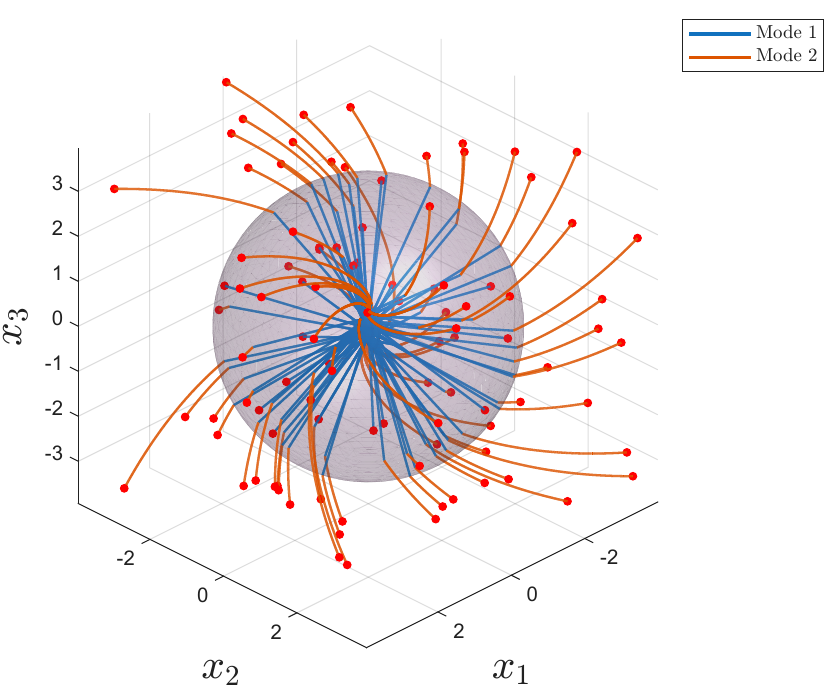}
      \caption{Trajectories of the joint state of the switching network system. Mode 1 and mode 2 are labeled with different colors, and each initial state is marked as red.}
      \label{fig: switching network traj colored}
    \end{figure}

    Figure \ref{fig: switching network traj colored} shows the trajectories and the switching surface of the joint state of this switching network system with parameters $\alpha_0 = 0$, $\alpha_1 = -0.1$, $\alpha_2 = 0.01$, and $R = 3.0$. 

\begin{figure}[H]
    \centering
    \small
    \setlength{\tabcolsep}{0pt}
    \begin{tabular}{@{}cc@{}}
    \includegraphics[width=.50\linewidth]{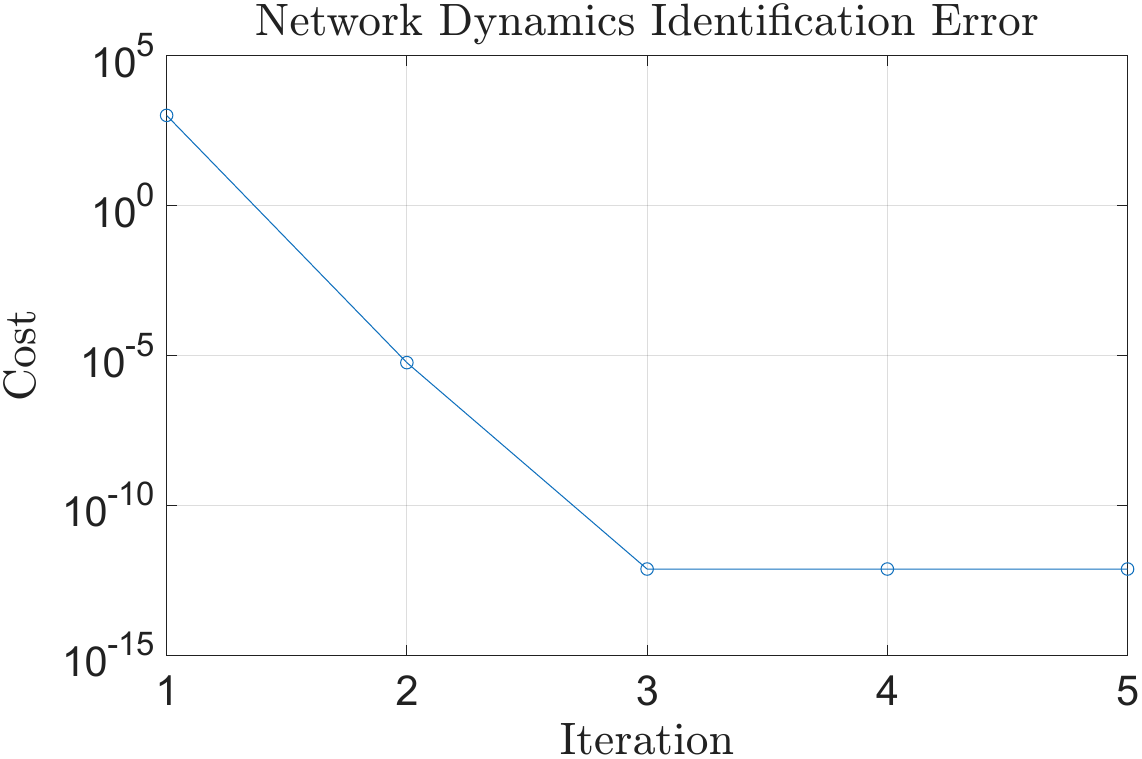} &
    \includegraphics[width=.50\linewidth]{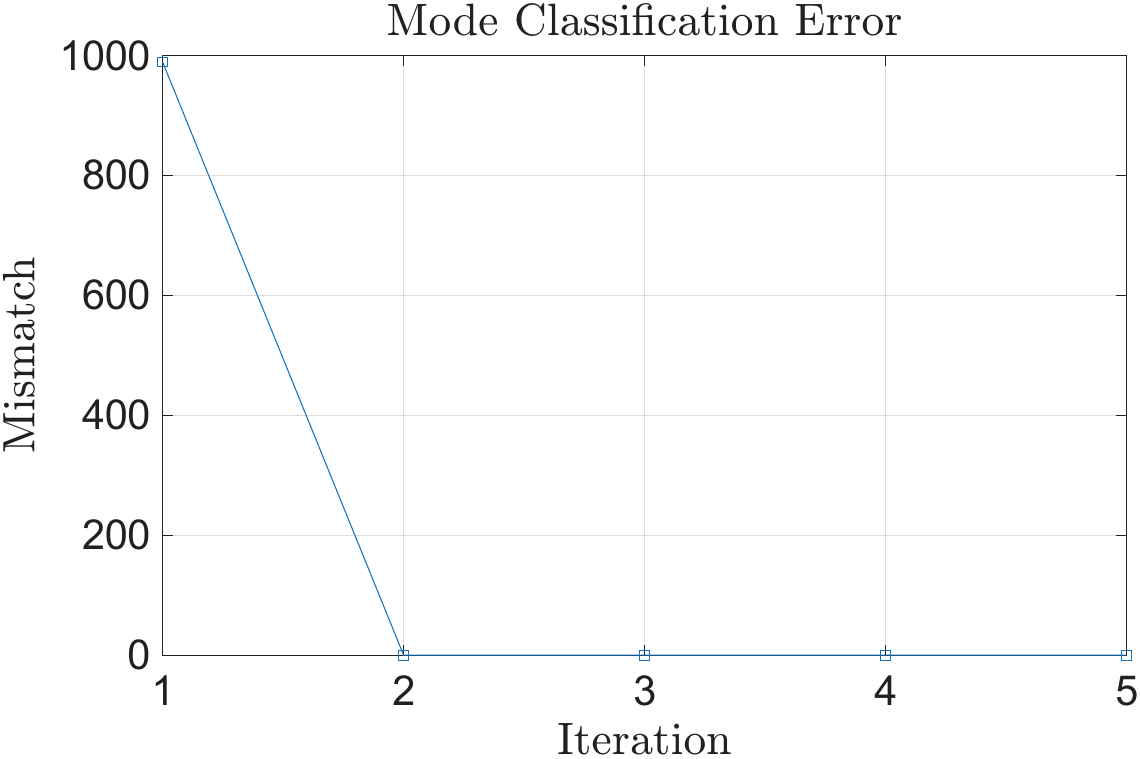} \\
    \includegraphics[width=.50\linewidth]{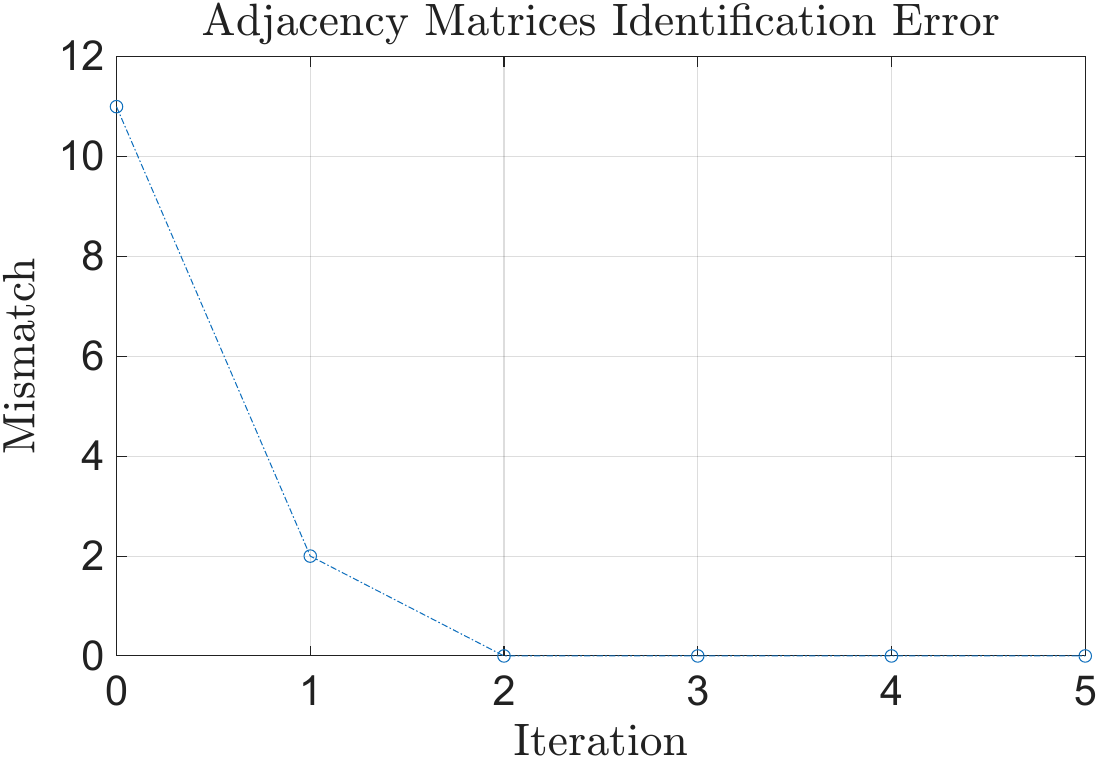} &
    \includegraphics[width=.50\linewidth]{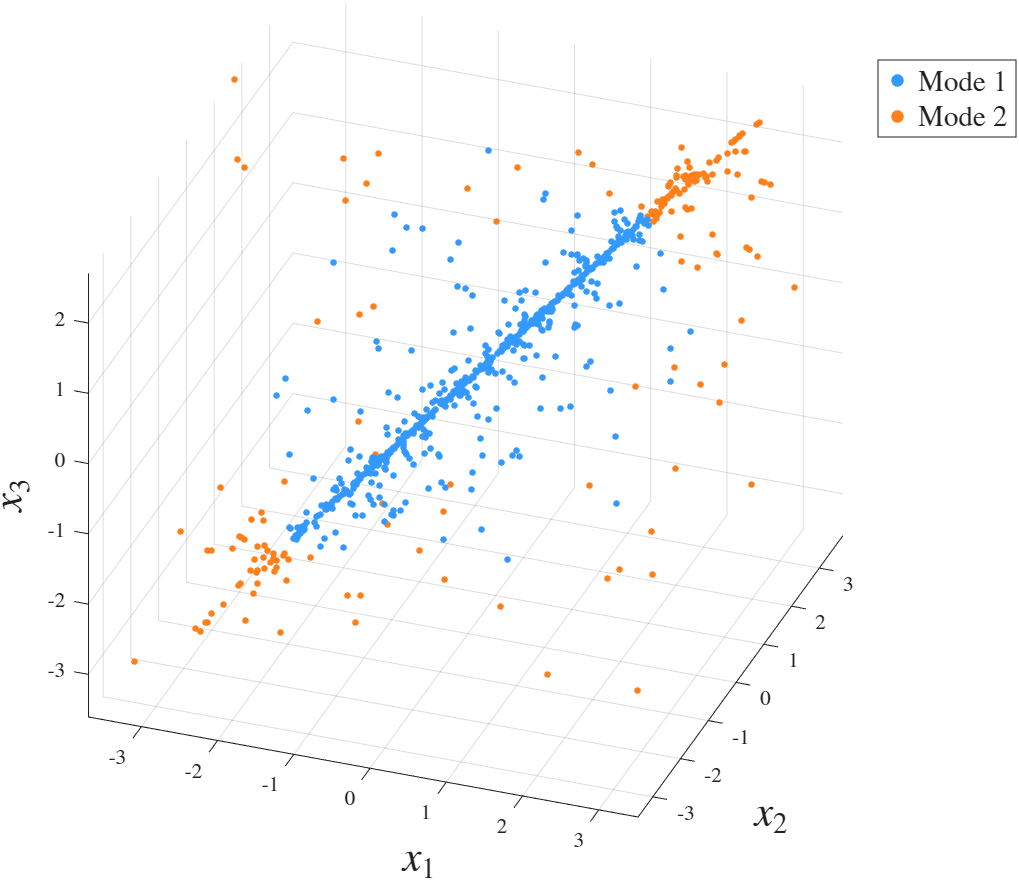}
    \end{tabular}
    \caption{Convergence and identification performance of the proposed bilevel convex optimization framework for switching network systems. Top left: Evolution of the network dynamics reconstruction error per iteration. Top right: Mode classification error showing convergence of mode assignments to the ground truth. Bottom left: Adjacency matrix identification error indicating recovery of the hidden network topology. Bottom right: 3D scatter plot of the sampled state data colored by identified modes.}
    \label{fig: identification errors}
    \end{figure}

\begin{figure}[H]
        \centering
        \includegraphics[width=0.85\linewidth]{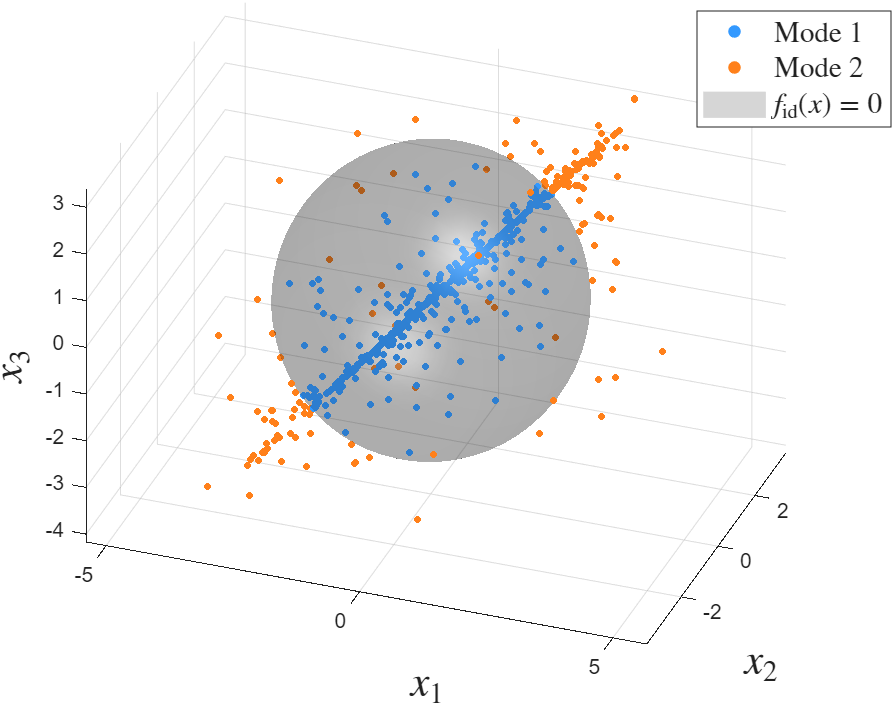}
        \caption{Identified switching surface $f_{\text{id}}(x)=0$ obtained via the soft-margin SVM approach in \cite{iwasaki2025learning}. The surface separates the two identified dynamic regimes in the state space.}
        \label{fig: identified_surface}
    \end{figure}

\begin{figure}[H]
      \centering
      \includegraphics[width=0.85\linewidth]{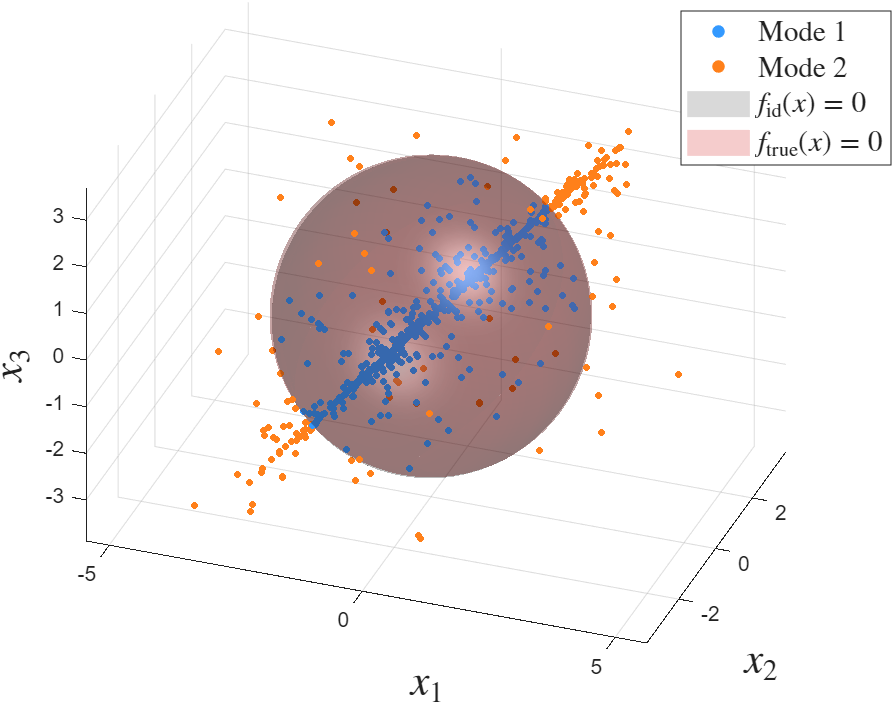}
      \caption{Comparison of recovered switching surface $f_{\text{id}}=0$ and ground truth $f_{\text{true}} = 0$. The identified switching surface well-separates the mode 1 and mode 2 samples, indicating it approximates the true switching rule well.}
      \label{fig: switching network comparison}
    \end{figure}

   Figure \ref{fig: identification errors} shows the convergence of bilevel convex optimization and the mismatch of mode and adjacency matrices from the ground truth and the scattering plot of the mode assignments to each sampled data point, highlighted in two different colors. 
    After identifying the mode assignments, we reconstructed the switching surface using the soft-margin SVM formulation described in our previous work \cite{iwasaki2025learning}. The classifier separates the mode-labeled samples in the continuous state space by learning a polynomial decision function $f_\text{id}(x)=0$, which defines the estimated switching boundary. Figure \ref{fig: identified_surface} visualizes the learned surface, while Figure \ref{fig: switching network comparison} compares it with the ground-truth switching law $f_{\text {true}}\left(x_1, x_2, x_3\right)=x_1^2+x_2^2+x_3^2-R^2$. The identified surface accurately separates the two modes and closely matches the analytic boundary, confirming that the proposed identification framework captures the underlying switching rule.
    
\section{Conclusion and future work}
This paper presented a bilevel convex optimization framework for identifying switching network systems from trajectory data. By alternating between mode assignment and subsystem identification, the proposed method recovers both the continuous dynamics and the underlying switching logic. Numerical experiments demonstrated that the identified mode dynamics, adjacency matrices, and switching surfaces closely match their ground-truth counterparts, confirming the effectiveness of the convex relaxation approach.

Future work will extend this framework in several directions. First, we aim to extend the proposed method to large-scale state spaces such as high-dimensional data and/or large number of mode structures. Another extension of this research could be to explore hypergraph-based formulations that capture higher-order couplings beyond pairwise interactions. On the learning side, to improve adaptability to larger class of network systems, developing model-free identification schemes is necessary. Finally, in the general hybrid system identification perspective, one needs to take account of systems where the switching logic depends not only on the instantaneous state or time, but also on its mode itself such as systems with hysteresis.

\addtolength{\textheight}{-12cm}   




\section*{Acknowledgements} Supported in part by NSF grant 2103026, and AFOSR
grants FA9550-32-1-0215 and FA9550-23-1-0400 (MURI).


{\footnotesize 
\balance
\bibliographystyle{IEEEtran}
\bibliography{bib/strings-abrv,bib/ieee-abrv,bib/references}

\begin{thebibliography}{10}
\providecommand{\url}[1]{#1}
\csname url@rmstyle\endcsname
\providecommand{\newblock}{\relax}
\providecommand{\bibinfo}[2]{#2}
\providecommand\BIBentrySTDinterwordspacing{\spaceskip=0pt\relax}
\providecommand\BIBentryALTinterwordstretchfactor{4}
\providecommand\BIBentryALTinterwordspacing{\spaceskip=\fontdimen2\font plus
\BIBentryALTinterwordstretchfactor\fontdimen3\font minus \fontdimen4\font\relax}
\providecommand\BIBforeignlanguage[2]{{%
\expandafter\ifx\csname l@#1\endcsname\relax
\typeout{** WARNING: IEEEtran.bst: No hyphenation pattern has been}%
\typeout{** loaded for the language `#1'. Using the pattern for}%
\typeout{** the default language instead.}%
\else
\language=\csname l@#1\endcsname
\fi
#2}}

\bibitem{iwasaki2025learning}
K.~Iwasaki, S.~Teng, A.~Bloch, and M.~Ghaffari, ``Learning hybrid dynamics via convex optimizations,'' \emph{arXiv preprint arXiv:2509.24157}, 2025.

\bibitem{mouyebe2025network}
M.~Mouyebe and A.~Bloch, ``Coupling induced stabilization of network dynamical systems and switching,'' \emph{arXiv preprint arXiv:2504.00403}, 2025.

\bibitem{Mirollo1990synchro}
R.~E. Mirollo and S.~H. Strogatz, ``Synchronization of pulse-coupled biological oscillators,'' \emph{SIAM Journal on Applied Mathematics}, vol.~50, no.~6, pp. 1645--1662, 1990.

\bibitem{Acebron2005Kuramoto}
J.~A. Acebr\'on, L.~L. Bonilla, C.~J. P\'erez~Vicente, F.~Ritort, and R.~Spigler, ``The kuramoto model: A simple paradigm for synchronization phenomena,'' \emph{Rev. Mod. Phys.}, vol.~77, pp. 137--185, Apr 2005.

\bibitem{Dorfler2013powergrids}
F.~Dörfler, M.~Chertkov, and F.~Bullo, ``Synchronization in complex oscillator networks and smart grids,'' \emph{Proceedings of the National Academy of Sciences}, vol. 110, no.~6, pp. 2005--2010, 2013.

\bibitem{Olfati2007consensus}
R.~Olfati-Saber, J.~A. Fax, and R.~M. Murray, ``Consensus and cooperation in networked multi-agent systems,'' \emph{Proceedings of the IEEE}, vol.~95, no.~1, pp. 215--233, 2007.

\bibitem{poli2020gnode}
M.~Poli, S.~Massaroli, J.~Park, A.~Yamashita, H.~Asama, and J.~Park, ``Graph neural ordinary differential equations,'' in \emph{Workshop on Deep Learning on Graphs: Methodologies and Applications (DLGMA'20)}, 2020.

\bibitem{moreau2005stability}
L.~Moreau, ``Stability of multiagent systems with time-dependent communication links,'' \emph{IEEE Transactions on Automatic Control}, vol.~50, no.~2, pp. 169--182, 2005.

\bibitem{yang2015network}
T.~Yang, J.~Wu, W.~Lu, and G.~Chen, ``Network synchronization with nonlinear dynamics and switching interactions,'' \emph{IEEE Transactions on Automatic Control}, vol.~61, no.~10, pp. 3103--3108, 2015.

\bibitem{kalofolias2016learn}
V.~Kalofolias, ``How to learn a graph from smooth signals,'' in \emph{Proceedings of the 19th International Conference on Artificial Intelligence and Statistics (AISTATS 2016)}, ser. Proceedings of Machine Learning Research, vol.~51.\hskip 1em plus 0.5em minus 0.4em\relax Journal of Machine Learning Research (JMLR), May 2016, pp. 920--929.

\bibitem{dong2019learn}
X.~Dong, D.~Thanou, M.~Rabbat, and P.~Frossard, ``Learning graphs from data: A signal representation perspective,'' \emph{IEEE Signal Processing Magazine}, vol.~36, no.~3, pp. 44--63, 2019.

\bibitem{timme2014revealing}
M.~Timme and J.~Casadiego, ``Revealing networks from dynamics: An introduction,'' \emph{Journal of Physics A: Mathematical and Theoretical}, vol.~47, no.~34, p. 343001, 2014.

\bibitem{ching2015reconstruct}
E.~S.~C. Ching, P.-Y. Lai, and C.~Y. Leung, ``Reconstructing weighted networks from dynamics,'' \emph{Phys. Rev. E}, vol.~91, p. 030801, Mar 2015.

\bibitem{CurtoFialkow1996}
R.~E. Curto and L.~A. Fialkow, \emph{Solution of the Truncated Complex Moment Problem for Flat Data}, ser. Memoirs of the American Mathematical Society.\hskip 1em plus 0.5em minus 0.4em\relax American Mathematical Society, 1996, vol. 568.

\bibitem{CurtoFialkow1998}
------, \emph{Flat Extensions of Positive Moment Matrices: Recursively Generated Relations}, ser. Memoirs of the American Mathematical Society.\hskip 1em plus 0.5em minus 0.4em\relax American Mathematical Society, 1998, vol. 648.

\bibitem{Lasserre2001Moments}
J.~B. Lasserre, ``Global optimization with polynomials and the problems of moments,'' \emph{SIAM Journal on Optimization}, vol.~11, no.~3, pp. 796--817, 2001.

\bibitem{Lasserre_2015}
------, \emph{An Introduction to Polynomial and Semi-Algebraic Optimization}, ser. Cambridge Texts in Applied Mathematics.\hskip 1em plus 0.5em minus 0.4em\relax Cambridge University Press, 2015.

\bibitem{Teng2025ACC}
S.~Teng, K.~Iwasaki, W.~Clark, X.~Yu, A.~Bloch, R.~Vasudevan, and M.~Ghaffari, ``A generalized metriplectic system via free energy and system identification via bilevel convex optimization,'' in \emph{2025 American Control Conference (ACC)}, 2025, pp. 1827--1833.

\end{thebibliography}
}

\end{document}